\documentclass[doublespacing]{elsart}

\usepackage{amsmath,amsfonts,amssymb}
\usepackage[dvips]{graphicx}
\usepackage[english,french]{babel}

\newtheorem{theorem}{Theorem}[section]

\newtheorem{e-proposition}[theorem]{Proposition}

\newtheorem{e-definition}[theorem]{Definition\rm}
\newtheorem{remark}{\it Remark\/}

\renewcommand{\theequation}{\arabic{equation}}
\def\og{\leavevmode\raise.3ex\hbox{$\scriptscriptstyle\langle\!\langle$~}}
\def\fg{\leavevmode\raise.3ex\hbox{~$\!\scriptscriptstyle\,\rangle\!\rangle$}}

\begin{document}
\begin{center}
\selectlanguage{english}
\title{ Navier-Stokes equations with periodic boundary conditions and pressure loss}
\end{center}
\begin{center}
\author[1]{Ch\'erif Amrouche},
\author[1,2]{Macaire Batchi},
\author[2]{Jean Batina}.
\address[1]{Laboratoire de Math\'ematiques Appliqu\'ees,
CNRS UMR 5142}
\address[2]{Laboratoire de Thermique Energ\'etique et
Proc\'ed\'es}
\address{Universit\'e de Pau et des Pays de l'Adour}
\address{Avenue de l'Universit\'e 64000 Pau, France}
\end{center}
\begin{center}
\textbf{Abstract}
\end{center}
\noindent We present in this note the existence and %
uniqueness results for the Stokes and Navier-Stokes equations %
which model the laminar flow of an incompressible fluid inside a %
two-dimensional channel of periodic sections.  The data of the %
pressure loss coefficient enables us to establish a relation on %
the pressure and to thus formulate an equivalent problem.
\selectlanguage{english}
\section{Introduction}
The problem which one proposes to study here is that modelling a
laminar flow inside a  two-dimensional plane channel with periodic
section. \noindent Let $\Omega $ be an open bounded connected
lipschtzian of $\Rset^2$\ (see figure hereafter), where
\qquad\ \ \ \ \ \ \ \ \ \ \ \ \ $\Gamma _{0}=\left\{ 0\right\}
\times \left] -1,1\right[ \text{ and }\Gamma _{1}=\left\{
1\right\} \times
\left] -1,1\right[.$\\
\noindent One defines the space
\begin{displaymath}
V=\left\{ \boldsymbol{v}\text{ }\mathbf{\in H}^{1}\left( \Omega \right) ;%
\text{div }\boldsymbol{v}\text{ }\mathbf{=}\text{ }0,\boldsymbol{v}\text{ }%
\mathbf{=0}\text{ on }\Gamma _{2},\text{ }\boldsymbol{v}_{\mid
_{\Gamma _{0}}}=\boldsymbol{v}_{\mid _{\Gamma _{1}}}\right\}
\end{displaymath}
\noindent and for $\lambda \in \mathbb{R} $ given, one considers
the problem
\begin{displaymath}
\left( \mathcal{S}\right) \left\{
\begin{array}{ccc}
\text{Find  }\boldsymbol{u}\in V,\text{ such that} \\
\forall \boldsymbol{v}\text{ }\mathbf{\in }\text{ }V,\text{ }%
\displaystyle\int_{\Omega }\text{ }\nabla
\boldsymbol{u}\mathbf{.}\nabla
\boldsymbol{v}\text{ }d\boldsymbol{x}=\lambda \displaystyle\int_{-1}^{+1}\text{ }%
v_{1}\left( 1,y\right) \text{ }dy.%
\end{array}%
\right.
\end{displaymath}
\begin{figure}[h]
\begin{center}
\includegraphics[width=6cm,height=4cm]{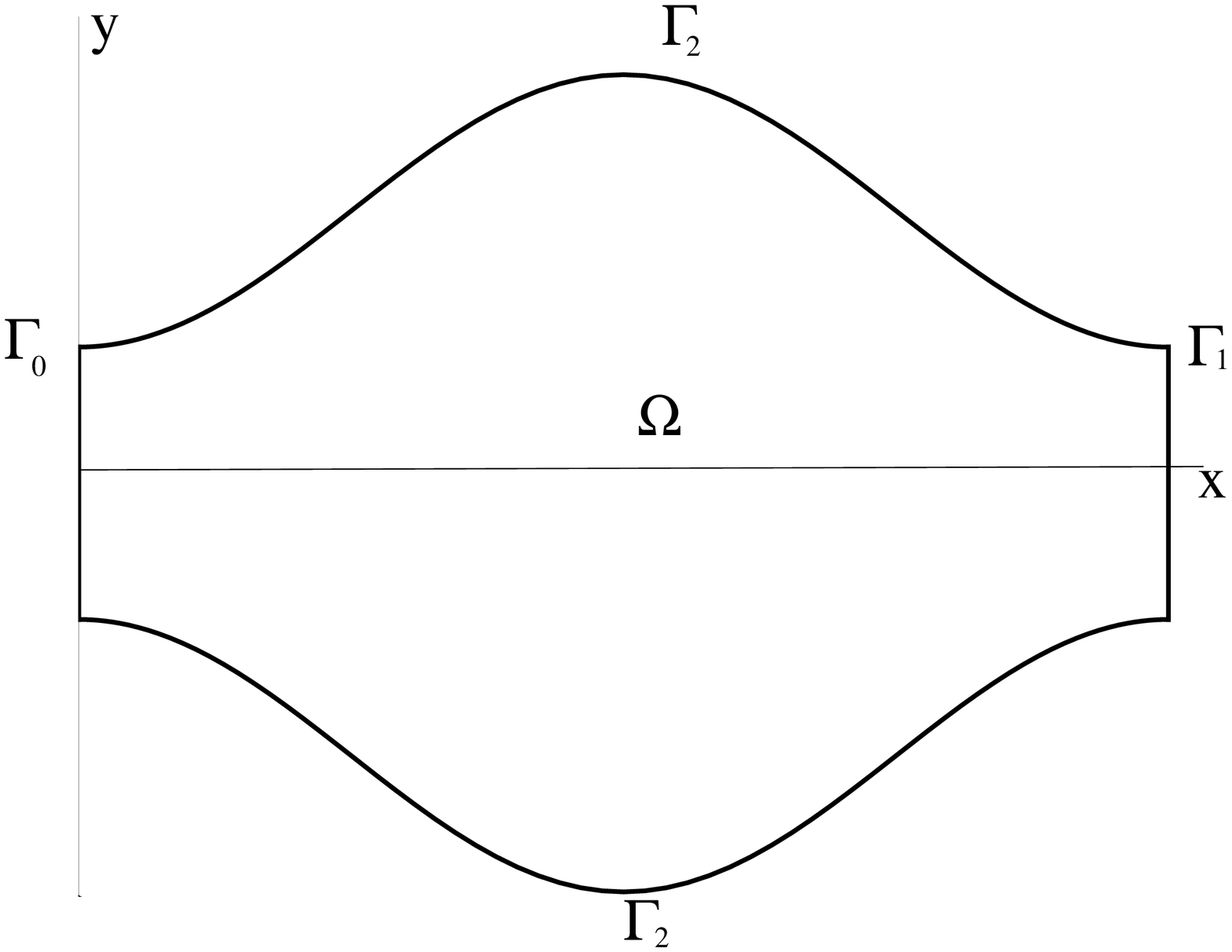}
\end{center}
\caption[]{Geometry of channel}
\end{figure}
\section{Resolution of the problem $\left( \mathcal{S}\right) $}
Initially one proposes to study the problem $\left( \mathcal{P}%
\right) .$ One has it\\
\begin{theorem}
\textit{Problem }$\left(
\mathcal{S}\right) $\textit{\ has an unique solution }$\boldsymbol{u%
}$ $\mathbf{\in }$ $V.$ \textit{Moreover, there is a constant
}$C\left( \Omega \right)
>0$ \textit{such that:}
\begin{equation}\label{eq}
\left\Vert \boldsymbol{u}\right\Vert _{\mathbf{H}^{1}\left( \Omega
\right) }\leq \lambda \text{C}\left( \Omega \right).
\end{equation}%
\end{theorem}
\noindent \textit{Proof: }Let us note initially that space $V$
provided the norm $H^{1} \left( \Omega \right) ^{2}$\quad%
being a closed subspace of %
$H^{1}\left( \Omega \right) ^{2}$ is thus an Hilbert space.
\noindent Let us set
\begin{displaymath}
a\left( \boldsymbol{u}\mathbf{,}\boldsymbol{v}\right)
=\int_{\Omega
}\nabla \boldsymbol{u}\mathbf{.}\nabla \boldsymbol{v}\text{ }d\boldsymbol{x}%
,\qquad l\left( \boldsymbol{v}\right) =\lambda
\int_{-1}^{+1}v_{1}\left( 1,y\right) \text{ }dy.
\end{displaymath}
\noindent It is clear, thanks to the Poincar\'e inequality, that
the bilinear continuous form is $V$-coercive.\ It is easy to also
see
that $l\in V^{\prime }.$\ One deduces from  Lax-Milgram Theorem  %
the existence and uniqueness of $\boldsymbol{u}$  solution of $%
\left( \mathcal{S}\right) .\ $\ Moreover,
\begin{displaymath}
\int_{\Omega }\text{ }\left\vert \nabla \boldsymbol{u}\right\vert
^{2}d\boldsymbol{x}\leq \lambda \sqrt{2}\left( \int_{-1}^{+1}\text{ }%
\left\vert u_{1}\left( 1,y\right) \right\vert ^{2}dy\right)
^{1/2},
\end{displaymath}%
\textit{i.e.}
\begin{displaymath}
\left\Vert \nabla \boldsymbol{u}\right\Vert _{L^{2}\left( \Omega
\right) }^{2} \leq \lambda \sqrt{2}\left\Vert
\boldsymbol{u}\right\Vert _{L^{2}\left( \Gamma \right) }\leq
\lambda \sqrt{2}\left\Vert \boldsymbol{u}\right\Vert
_{H^{1/2}\left( \Gamma \right) }
\end{displaymath}
\begin{displaymath}
\leq \lambda C_{1}\left( \Omega \right) \left\Vert
\boldsymbol{u}\right\Vert _{H^{1}\left( \Omega \right) }
\end{displaymath}
Thus there is the estimate (\ref{eq}).$\square \medskip $
\\
\noindent We now will give an interpretation of the problem
$\left( \mathcal{S}\right) .$\ One introduces the space
\begin{displaymath}
\mathcal{V}=\left\{ \boldsymbol{v}\text{ }\mathbf{\in }\text{ }\mathcal{D}%
\left( \Omega \right) ^{2};\quad\text{div }\boldsymbol{v}\text{
}\mathbf{=}\text{ }0\right\} .
\end{displaymath}
\noindent Let $\boldsymbol{u}$ be the solution of $\left(
\mathcal{S}\right) .$ Then, for all $\boldsymbol{v}$ $\mathbf{\in
}$ $\mathcal{V},$ one has :
\begin{displaymath}
\left\langle -\Delta \boldsymbol{u}\mathbf{,}\text{ }
\boldsymbol{v}\right\rangle _{\mathcal{D}^{\prime }\left( \Omega
\right) \times \mathcal{D}\left( \Omega \right) }=0.
\end{displaymath}
\noindent So that thanks to De Rham Theorem, there exists $%
p\in \mathcal{D}^{\prime }\left( \Omega \right) $ such that
\begin{equation}
-\Delta \boldsymbol{u}+\mathbf{\nabla }p=0\text{ in }\Omega .
\end{equation}
\noindent Moreover, since $\ \mathbf{\nabla }p\in H^{-1}\left(
\Omega \right) ^{2},$ it is known that there exists $q\in $
$L^{2}\left( \Omega \right) $ such that (see  $\left[ 1\right] $)
\begin{equation}
\mathbf{\nabla }q=\mathbf{\nabla }p\quad\text{ in }\Omega .
\end{equation}
\noindent The open $\Omega $ being connected, there exists
$C\in \mathbb{R} $ such that $p=q+C,$ what means that $p\in $ $L^{2}\left( \Omega \right) $%
. Let us recall that (see $\left[ 1\right] $)
\begin{displaymath}
\underset{K\text{ }\in \text{ }%
\mathbb{R}
}{\inf }\left\Vert p+K\right\Vert _{L^{2}\left( \Omega \right)
}\leq C\left\Vert \mathbf{\nabla }p\right\Vert _{H^{-1}\left(
\Omega \right) ^{2}}.
\end{displaymath}
\noindent One deduces from the estimate (1) and from (2) that
\begin{eqnarray*}
\underset{K\text{ }\in \text{ }%
\mathbb{R}
}{\inf }\left\Vert p+K\right\Vert _{L^{2}\left( \Omega \right) }
&\leq &C\left\Vert \Delta \boldsymbol{u}\right\Vert _{H^{-1}\left(
\Omega \right) ^{2}}\leq C\left\Vert \boldsymbol{u}\right\Vert
_{H^{1}\left( \Omega \right) ^{2}} \leq \lambda \text{C}\left(
\Omega \right) .
\end{eqnarray*}
\noindent Since $\boldsymbol{u}$ $\mathbf{\in }$ $H^{1}\left(
\Omega \right) ^{2}$ and $\ \mathbf{0}=-\Delta
\boldsymbol{u}+\mathbf{\nabla }p\in L^{2}\left( \Omega \right)
^{2},$ it is shown that $\ -\dfrac{\partial
\boldsymbol{u}}{\partial \boldsymbol{n}}+p\boldsymbol{n}$ $\mathbf{\in }$ $%
H^{-1/2}\left( \Gamma \right) ^{2}$ and one has the Green formula: for all $%
\boldsymbol{v}$ $\mathbf{\in }$ $V$
\begin{equation}
\int_{\Omega }\left( -\triangle \boldsymbol{u}+\mathbf{\nabla }%
p\right) .\boldsymbol{v}\text{ }d\boldsymbol{x}=\int_{\Omega
}\nabla\boldsymbol{u}\mathbf{.}\nabla \boldsymbol{v}\text{ }d\boldsymbol{x}%
+\left\langle -\dfrac{\partial \boldsymbol{u}}{\partial \mathbf{n}}+p\boldsymbol{%
n,\text{ }}\boldsymbol{v}\right\rangle ,
\end{equation}
\noindent where the bracket represents the duality product
$H^{-1/2}\left(
\Gamma \right) \times H^{1/2}\left( \Gamma \right) .$ Moreover, as $%
p\in L^{2}\left( \Omega \right) $ and $\triangle p=0$ in $\Omega ,$ one has $%
p\in H^{-1/2}\left( \Gamma \right) .$ Consequently, one has therefore $\dfrac{%
\partial \boldsymbol{u}}{\partial \boldsymbol{n}}\in H^{-1/2}\left( \Gamma
\right) ^{2}.\smallskip $ \noindent The function $\boldsymbol{u}$
being solution of $\left( \mathcal{S}\right) , $ for all
$\boldsymbol{v}$ $\mathbf{\in }$ $V$\text{ }one has:%
\begin{equation}
\left\langle \dfrac{\partial \boldsymbol{u}}{\partial \boldsymbol{n}}-p\boldsymbol{n,%
}\boldsymbol{v}\right\rangle =\lambda \int_{-1}^{+1}v_{1}\left(
1,y\right) \text{ }dy,
\end{equation}
\noindent\textit{i.e.}
\begin{equation}
\left\langle \dfrac{\partial \boldsymbol{u}}{\partial x}-p\mathbf{e}_{1}%
\mathbf{,}\boldsymbol{v}\right\rangle _{\Gamma _{1}}+\left\langle -\dfrac{%
\partial \boldsymbol{u}}{\partial x}+p\mathbf{e}_{1}\mathbf{,}\boldsymbol{v}%
\right\rangle _{\Gamma _{0}}=\left\langle \lambda \mathbf{e}_{1,}\boldsymbol{v}%
\right\rangle _{\Gamma _{1}}.
\end{equation}%
\noindent where $\mathbf{e}_{1}=\left(1,0\right)$.
\smallskip\\
\noindent i) Let $\mu \in H_{00}^{1/2}\left( \Gamma _{1}\right) $ $%
\smallskip $and let us set
\begin{equation*}
\mu _{2}=\left\{
\begin{array}{lll}
\mu  & \text{on} & \Gamma _{0}\cup \Gamma _{1} \\
0 & \text{on} & \Gamma _{2}%
\end{array}%
\right.\quad \text{ and}\quad\ \boldsymbol{\mu =}\left(
\begin{array}{c}
0 \\
\mu _{2}%
\end{array}%
\right)
\end{equation*}
\noindent where (see  $\left[ 2\right] $)
\begin{equation*}
H_{00}^{1/2}\left( \Gamma _{1}\right) =\left\{ \varphi \in \mathbf{L}%
^{2}(\Gamma _{1})\mathbf{;}\text{ }\mathbf{\exists }\text{ }\boldsymbol{v}%
\text{ }\mathbf{\in H}^{1}\mathbf{(}\Omega \mathbf{),}\text{ with }%
\boldsymbol{v}\mathbf{\mid }_{\Gamma _{2}}\mathbf{=0,}\text{ }\boldsymbol{v}%
\mathbf{\mid }_{\Gamma _{0}\cup \Gamma _{1}}\text{ }\mathbf{=}\text{ }%
\varphi \right\}.
\end{equation*}
\noindent It is checked easily that
\\
\begin{equation*}
\boldsymbol{\mu \in }\text{ }H^{1/2}\left( \Gamma \right)
^{2}\quad\text{ and }\quad %
\int_{\Gamma }\boldsymbol{\mu .n}\text{
}d\sigma =0.
\end{equation*}
\noindent So that there exists $\boldsymbol{v}$ $\mathbf{\in }$ $%
H^{1}\left( \Omega \right) ^{2}$ satisfying (see $%
\left[ 3\right] $)
\begin{equation*}
\begin{array}{lllllll}
\text{div }\boldsymbol{v}\text{ }\mathbf{=}\text{ }0 &
\quad\text{in} & \quad\Omega \quad & \text{ and }\quad &
\boldsymbol{v}=\boldsymbol{\mu } & \quad\text{on} & \quad\Gamma.
\end{array}%
\end{equation*}
\noindent In particular $\boldsymbol{v}$ $\mathbf{\in }$ $V$ and
according to (6), one has
\begin{equation*}
\left\langle \dfrac{\partial u_{2}}{\partial x}\mathbf{,}\text{
}\mu
\right\rangle _{\Gamma _{1}}=\left\langle \dfrac{\partial u_{2}}{\partial x}%
\mathbf{,}\text{ }\mu \right\rangle _{\Gamma _{0}},
\end{equation*}
\noindent which means that
\begin{equation}
\dfrac{\partial u_{2}}{\partial x}\mathbf{\mid }_{\Gamma _{1}}=\dfrac{%
\partial u_{2}}{\partial x}\mid _{\Gamma _{0}}.
\end{equation}
\noindent One deduces now from (6) that for all $\boldsymbol{v}$ $%
\mathbf{\in }$ $V,$
\begin{equation}
\left\langle \dfrac{\partial u_{1}}{\partial x}-p\mathbf{,}\text{ }%
v_{1}\right\rangle _{\Gamma _{1}}+\left\langle -\dfrac{\partial u_{1}}{%
\partial x}+p\mathbf{,}\text{ \ }v_{1}\right\rangle _{\Gamma
_{0}}=\left\langle \lambda \mathbf{,}\text{ }v_{1}\right\rangle
_{\Gamma _{1}}.
\end{equation}
\noindent But, div $\boldsymbol{u}=0$ and \ $u_{2}\mathbf{\mid
}_{\Gamma _{1}}=u_{2}\mid _{\Gamma _{0}},$ one thus has

\begin{equation}
\dfrac{\partial u_{2}}{\partial y}\mathbf{\mid }_{\Gamma _{1}}=\dfrac{%
\partial u_{2}}{\partial y}\mid _{\Gamma _{0}}\quad\text{and}\quad%
\dfrac{\partial u_{1}}{\partial x}\mathbf{\mid }_{\Gamma _{1}}=\dfrac{%
\partial u_{1}}{\partial x}\mid _{\Gamma _{0}}.
\end{equation}
\noindent Consequently, thanks to (8) one has:
\begin{equation}
\left\langle -p\mathbf{,}\text{ }v_{1}\right\rangle _{\Gamma
_{1}}+\left\langle p\mathbf{,}\text{ }v_{1}\right\rangle _{\Gamma
_{0}}=\left\langle \lambda \mathbf{,}\text{ }v_{1}\right\rangle
_{\Gamma _{1}}
\end{equation}
\noindent ii) Let $\nu \in H_{00}^{1/2}\left( \Gamma _{1}\right) $ and let us set%
\begin{equation*}
\nu _{1}=\left\{
\begin{array}{lll}
\nu & \text{on} & \Gamma _{0}\cup \Gamma _{1} \\
0 & \text{on} & \Gamma _{2}%
\end{array}%
\right.\quad \text{ and }\quad\ \boldsymbol{\nu}=\left(
\begin{array}{c}
\nu _{1} \\
0%
\end{array}%
\right) .
\end{equation*}%
One easily checks that
\begin{equation*}
\boldsymbol{\nu \in }\text{ }H^{1/2}\left( \Gamma \right)
^{2}\quad\text{ and }\quad %
\int_{\Gamma }\boldsymbol{\nu .n}\text{
}d\sigma =0.
\end{equation*}
\noindent So that there exists $\boldsymbol{v}$ $\mathbf{\in }$ $%
H^{1}\left( \Omega \right) ^{2}$ satisfying%
\begin{equation*}
\begin{array}{lllllll}
\text{div }\boldsymbol{v}\text{ }\mathbf{=}\text{ }0 &
\quad\text{in} & \quad\Omega
& \quad \text{ and }\quad & \boldsymbol{v}=\boldsymbol{\nu }\quad & \text{on} & \quad\Gamma .%
\end{array}%
\end{equation*}
\noindent In particular $\boldsymbol{v}$ $\mathbf{\in }$ $V$ and
according to (14), one has
\begin{equation*}
\left\langle -p\mathbf{,}\text{ }\nu \right\rangle _{\Gamma
_{1}}+\left\langle p\mathbf{,}\text{ }\nu \right\rangle _{\Gamma
_{0}}=\left\langle \lambda \mathbf{,}\text{ }\nu \right\rangle
_{\Gamma _{1}}
\end{equation*}%
\textit{i.e.}
\begin{equation}
p_{\mathbf{\mid }\Gamma _{1}}=p_{\mathbf{\mid }\Gamma
_{0}}-\lambda
\end{equation}
\smallskip\\
\noindent where the equality takes place with the  $H^{1/2} $\quad
sense. \noindent In short, if $\boldsymbol{u}$ $\mathbf{\in }$ $%
H^{1}\left( \Omega \right) ^{2}$ is solution of $\left(
\mathcal{S}\right) ,$ then there exists $p\in L^{2}\left( \Omega
\right) ,$ unique with an  additive constant such that:
\begin{equation}
\qquad\ \ \ \ -\Delta \boldsymbol{u}+\mathbf{\nabla
}p=\mathbf{0}\quad\text{in}\quad\Omega ,%
\end{equation}
\begin{equation}
\qquad \ \ \ \ \text{div }\boldsymbol{u}=0\ \ \ \ \ \ \ \  \quad\text{in}  \quad\Omega ,%
\end{equation}
\begin{equation}
\qquad\ \ \ \ \ \ \qquad\ \ \ \ \boldsymbol{u}\text{ }\mathbf{=0\ \ \  } \quad\text{on\ } \quad%
\Gamma
_{2,\quad}\\%
\boldsymbol{u}\mathbf{\mid }_{\Gamma _{1}}=\boldsymbol{u}\mid
_{\Gamma _{0}},
\end{equation}
\begin{equation}
\dfrac{\partial \boldsymbol{u}}{\partial x}\mathbf{\mid }_{\Gamma _{1}}=%
\dfrac{\partial \boldsymbol{u}}{\partial x}\mid _{\Gamma _{0}},
\end{equation}
\begin{equation}
p_{\mathbf{\mid }\Gamma _{1}}=p_{\mathbf{\mid }\Gamma
_{0}}-\lambda .
\end{equation}
\\
\noindent It is clear that if $\left( \boldsymbol{u},p\right) \in
H^{1}\left( \Omega \right) ^{2}\times L^{2}\left( \Omega \right) $
checks (12)-(16),\ %
then  $\boldsymbol{u}$ is solution of $\left( \mathcal{S%
}\right)$. \noindent Thus it\\
\medskip
\begin{theorem}
\textit{The problem (12)-(16)
has an unique solution  }$\left( \boldsymbol{u},p\right) \in $ $%
\mathit{H}^{1}\left( \Omega \right) ^{2}\times
\mathit{L}^{2}\left( \Omega \right) $\textit{\, }\textit{\ \ up to
an additive constant for $p$.\ Moreover,
}$\boldsymbol{u}$\textit{\ verifies }$\left( \mathcal{S}\right)
$\textit{\ and }
\begin{equation*}
\left\Vert \boldsymbol{u}\right\Vert _{\mathbf{H}^{1}\left( \Omega
\right)
}+\left\Vert p\right\Vert _{L^{2}\left( \Omega \right) /\text{ }%
\mathbb{R}
}\leq \lambda C\left( \Omega \right) .\text{ }\square
\smallskip\\
\end{equation*}
\end{theorem}
\begin{remark}
The pressure verifies the relation (16),
which means that $p$ satisfies the relation of Patankar et \textit{al}.\text{}$\left[ 5%
\right] .\smallskip $
\end{remark}
\section{ Navier-Stokes Equations}
One takes again the assumptions of the Stokes problem given above.
For $\lambda \in
\mathbb{R}
$ given, the one considers the following problem
\begin{equation*}
\left( \mathcal{NS}\right) \left\{
\begin{array}{c}
\text{Find $\boldsymbol{u}$}\in V,\text{ such that} \\
\forall \boldsymbol{v}\text{ }\mathbf{\in }\text{ }V,\text{ }%
\displaystyle\int_{\Omega }\text{ }\nabla
\boldsymbol{u}\mathbf{.}\nabla \boldsymbol{v}\text{
}d\boldsymbol{x}+b\left( \boldsymbol{u}\mathbf{,\ }\boldsymbol{u}\mathbf{,\ }\boldsymbol{v}%
\right)=\lambda \displaystyle\int_{-1}^{+1}\text{ }v_{1}\left(
1,y\right)
\text{ }dy%
\end{array}%
\right.
\end{equation*}
with
\begin{equation*}
b\left( \boldsymbol{u}\mathbf{,\ }\boldsymbol{v}\mathbf{,\ }\boldsymbol{w}%
\right) =\int_{\Omega }\text{ }\left( \boldsymbol{u}\mathbf{.\nabla }%
\right) \boldsymbol{v}\mathbf{.}\boldsymbol{w}\text{
}d\boldsymbol{x}
\end{equation*}
\smallskip\\
With an aim of establishing the existence of the solutions of the problem  $%
\left( \mathcal{NS}%
\right) ,$ one uses the Brouwer fixed point theorem (see $\left[
4\right] $,  $\left[ 6\right] $). One will show it
\smallskip
\begin{theorem}
\textit{The problem} $\left(
\mathcal{NS}\right) $ \textit{has at least a solution} $%
\boldsymbol{u}$ $\mathbf{\in }$ $V.$ \textit{Moreover,}
$\boldsymbol{u}$ \textit{checks the estimate (1).\smallskip }
\end{theorem}
\noindent \textit{Proof: }To show the existence of $%
\boldsymbol{u}$, one constructs the approximate solutions of the problem $\left( \mathcal{NS}\right)   $%
\text{ }by the Galerkin method and then thanks to the arguments of
compactness, one makes a passage to the limit.
\smallskip
\\
\noindent \textit{i) For each fixed integer }$m\geq 1,$
\textit{one defines an
 approximate solution }$\boldsymbol{u}_{m}$\textit{\ of }$\left( \mathcal{NS}%
\right) $\textit{\ by}
\begin{equation}
\begin{array}{c}
\boldsymbol{u}_{m}=\displaystyle\sum_{i=1}^{m}g_{im}\boldsymbol{w}_{i}\text{,
\quad\text{with} }\quad g_{im}\in
\mathbb{R}\\
\left( \left( \boldsymbol{u}_{m},\boldsymbol{w}_{i}\right)
\right) +b\left( \boldsymbol{u}_{m}\mathbf{,\  }\boldsymbol{u}_{m}\mathbf{,\  }%
\boldsymbol{w}_{i}\right) =\left\langle \lambda \boldsymbol{n,\ }\boldsymbol{w}%
_{i}\right\rangle _{\Gamma _{1}},i=1,...,m
\end{array}
\end{equation}%
\smallskip\\
\noindent where $V_{m}=$ $\left\langle \boldsymbol{w}_{1},...,\boldsymbol{w}%
_{m}\right\rangle $ vector spaces spanned by the vectors $%
\boldsymbol{w}_{1},...,\boldsymbol{w}_{m}$ and $\left\{ \boldsymbol{w}%
_{i}\right\} $ is an Hilbertian basis of $V$ which is
separable.\text{ }\smallskip \noindent Let us note that (17) is
equivalent to:
\begin{equation}
\forall \boldsymbol{v}\text{ }\mathbf{\in }\text{ }V_{m},\text{
}\left(
\left( \boldsymbol{u}_{m},\boldsymbol{v}\right) \right) +b\left( \boldsymbol{%
u}_{m}\mathbf{,\ }\boldsymbol{u}_{m},\  \boldsymbol{v}\right) =\text{ }%
\lambda \int_{-1}^{+1}\text{ }v_{1}\left( 1,y\right) \text{ }dy.
\end{equation}
\noindent With an aim to establish the existence of the solutions
of the problem $\boldsymbol{u}_{m},$ the operator as follows is
considered
\begin{equation*}
\begin{array}{llll}
\ \text{P}_{m}: & V_{m} & \rightarrow & V_{m} \\
& \boldsymbol{u} & \longmapsto & \text{P}_{m}\left( \boldsymbol{u}\right)%
\end{array}%
\end{equation*}
\noindent defined by \ \ \ \ \ \ $\ \ \ \ \ $%
\begin{equation*}
\forall \boldsymbol{u,v}\text{ }\mathbf{\in }\text{ }V_{m},\quad
\left( \left( \text{P}_{m}\left( \boldsymbol{u}\right)
,\boldsymbol{v}\right) \right) =\left( \left( \boldsymbol{u},
\boldsymbol{v}\right) \right) +b\left( \boldsymbol{u}\mathbf{, \
}\boldsymbol{u}\mathbf{,\ }\boldsymbol{v}\right) -\lambda
\int_{-1}^{+1}v_{1}\left( 1,y\right) \text{ }dy.
\end{equation*}
\\
\noindent Let us note initially that $\text{P}_{m}$ is continuous
and
$\ \ \ $%
\begin{equation*}
\forall \boldsymbol{u}\mathbf{ }\text{ }\mathbf{\in }\text{ }V,%
\text{ \quad }b\left( \boldsymbol{u}\mathbf{,\ }\boldsymbol{u}\mathbf{,\ }%
\boldsymbol{u}\right) =0.
\end{equation*}
\noindent Indeed, thanks to the Green formula, one has%
\begin{equation*}
\begin{array}{ll}
b\left( \boldsymbol{u}\mathbf{,\ }\boldsymbol{u}\mathbf{, \
}\boldsymbol{u}\right) & = -\frac{1}{2}\displaystyle\int_{\Omega
}\left\vert \boldsymbol{u}\right\vert
^{2}\text{div }\boldsymbol{u}\text{ }d\boldsymbol{x}+\frac{1}{2}%
\displaystyle\int_{\Gamma }\left(
\boldsymbol{u}\mathbf{\boldsymbol{.n}}\right)
\text{ }\left\vert \boldsymbol{u}\right\vert ^{2}d\sigma\text{}=0.%
\end{array}%
\end{equation*}
\noindent But, div $\boldsymbol{u}$ $=0$ in $\Omega $ and
\begin{equation*}
\begin{array}{lll}
\displaystyle\int_{\Gamma }$ $\left( \boldsymbol{u}\mathbf{\boldsymbol{.n}}\right) $ %
$\left\vert \boldsymbol{u}\right\vert ^{2}d\sigma
=\displaystyle\int_{\Gamma _{0}}$
$\left( \boldsymbol{u}\mathbf{\boldsymbol{.n}}\right) $ $\left\vert \boldsymbol{u%
}\right\vert ^{2}d\sigma +\displaystyle\int_{\Gamma _{1}}$ $\left( \boldsymbol{u}%
\mathbf{\boldsymbol{.n}}\right) $ $\left\vert
\boldsymbol{u}\right\vert ^{2}d\sigma.
\end{array}
\end{equation*}

\noindent since %
the external normal to $\Gamma _{0}$ is opposed to that of $\Gamma
_{1}$ and $\boldsymbol{u}$ $\mathbf{\in }$ $V$.
$\ \ \ $%
\\
\noindent Thanks to Brouwer Theorem, there exists
$\boldsymbol{u}_{m} $ satisfying (18) and

\begin{equation}
\left\Vert \boldsymbol{u}_{m}\right\Vert
_{\mathbf{H}^{1}\left(\Omega \right)} \leq \lambda C\left( \Omega
\right) .
\end{equation}
\smallskip\\
\noindent \textit{ii)} We can extract a subsequence
$\boldsymbol{u}_\nu$  such that
\smallskip\\
\begin{displaymath}
\boldsymbol{u}_\nu
\rightharpoonup\boldsymbol{u}\quad\text{weakly}\quad\text{in}\quad
V,
\end{displaymath}
\noindent and thanks to the compact imbedding of $V$ \text{in }  $%
L^{2}\left( \Omega \right) ^{2},$  we obtain
\begin{equation*}
\text{ }\forall \boldsymbol{v}\text{ }\mathbf{\in }\text{
}V\mathbf{,}\text{ }\left( \left(
\boldsymbol{u},\boldsymbol{v}\right) \right) +b\left(
\boldsymbol{u}\mathbf{,\ }\boldsymbol{u}\mathbf{,\
}\boldsymbol{v}\right) =\text{
}\lambda \int_{-1}^{+1}\text{ }v_{1}\left( 1,y\right) \text{ }dy\mathbf{.%
}
\end{equation*}
\noindent As for the Stokes problem, one shows the existence of $%
p\in L^{2}\left( \Omega \right) ,$ unique except for an  additive constant, such that%
\begin{equation*}
\left\{
\begin{array}{lll}
-\Delta \boldsymbol{u}+\left( \boldsymbol{u}\mathbf{.\nabla
}\right) \boldsymbol{u}\mathbf{\mathbf{+}\nabla }p=\mathbf{0} &
\quad\text{in}\quad & \Omega ,
\\
\text{div }\boldsymbol{u}=0 & \quad \text{in}\quad & \Omega , \\
\boldsymbol{u}\text{ }\mathbf{=0} & \quad \text{on}\quad & \Gamma _{2}, \\
\boldsymbol{u}\mathbf{\mid }_{\Gamma _{1}}=\boldsymbol{u}\mid
_{\Gamma _{0}}. &  &
\end{array}%
\right.
\end{equation*}
\noindent It is checked finally that
\begin{equation*}
\qquad\ \ \dfrac{\partial \boldsymbol{u}}{\partial x}\mathbf{\mid }_{\Gamma _{1}}=%
\dfrac{\partial \boldsymbol{u}}{\partial x}\mid _{\Gamma
_{0}}\text{ ,\ \ \ \ \ \ \ \ \ \ }
\end{equation*}%
\begin{equation*}
p_{\mathbf{\mid }\Gamma _{1}}=p_{\mathbf{\mid }\Gamma
_{0}}-\lambda .
\end{equation*}
\smallskip\\
\begin{remark}
i) Theorem 3.1 of problem $\left( \mathcal{NS}\right) $ takes
place in three dimension.
\\
ii) One can show that the solution $\left( \boldsymbol{u},p\right) \text{belongs to} $ $%
\mathit{H}^{2}\left( \Omega \right) ^{2}\times
\mathit{H}^{1}\left( \Omega \right) $.
\end{remark}


\end{document}